\documentclass[11pt]{article}

\usepackage{amssymb}
\input psfig.sty

\newcommand{\maths}[1]{{\mathbb #1}}  

\newtheorem{defi}{Definition}[section]
\newtheorem{prop}[defi]{Proposition}
\newtheorem{theo}[defi]{Theorem}
\newtheorem{conj}[defi]{Conjecture}
\newtheorem{lemm}[defi]{Lemma}
\newtheorem{coro}[defi]{Corollary}
\newtheorem{rema}[defi]{Remark}
\newtheorem{exem}[defi]{Example}

\newcommand{\bdefi}{\begin{defi}}
\newcommand{\edefi}{\end{defi}}
\newcommand{\bprop}{\begin{prop}}
\newcommand{\eprop}{\end{prop}}
\newcommand{\btheo}{\begin{theo}}
\newcommand{\etheo}{\end{theo}}
\newcommand{\blemm}{\begin{lemm}}
\newcommand{\brema}{\begin{rema}}
\newcommand{\erema}{\end{rema}}
\newcommand{\bexer}{\begin{exem}}
\newcommand{\eexer}{\end{exem}}
\newcommand{\bconj}{\begin{conj}}
\newcommand{\econj}{\end{conj}}
\newcommand{\elemm}{\end{lemm}}
\newcommand{\bcoro}{\begin{coro}}
\newcommand{\ecoro}{\end{coro}}
\newcommand{\dem}{\noindent{\bf Proof. }}

\renewcommand{\O}{{\cal O}}

\newcommand{\ra}{\rightarrow}

\newcommand{\RR}{\maths{R}}
\newcommand{\NN}{\maths{N}}

\newcommand{\QQ}{\maths{Q}}

\newcommand{\HH}{\maths{H}}

\newcommand{\ZZ}{\maths{Z}}

\newcounter{fig}

\newcommand{\PSLZ}{\mbox{${\rm{ PSL}}_{2}(\ZZ)$}} 
 
\newcommand{\PSLO}{\mbox{${\rm{ PSL}}_{2}(\O)$}}

\newcommand{\eop}[1]{{\hfill\fbox{\bf #1}}}

\title{Counting horoballs and rational geodesics} 

\author{Sa'ar Hersonsky \and Fr\'ed\'eric Paulin}

\date{\today}

\begin{document}
\maketitle

\begin{abstract}
\noindent 
Let $M$ be a geometrically finite pinched negatively curved Riemannian
manifold with at least one cusp.  We study the asymptotics of the
number of geodesics in $M$ starting from and returning to a given
cusp, and of the number of horoballs at parabolic fixed points in the
universal cover of $M$.  In the appendix, due to K.~Belabas, the case
of SL$(2,\ZZ)$ and of Bianchi groups is developed. \footnote{
rational geodesic, negative curvature, cusp, horoball.}
\end{abstract}

\section{Introduction}
\label{sec:Intro}

Let $M$ be a non elementary geometrically finite pinched negatively
curved Riemannian manifold with at least one cusp.  A geodesic line
starting from a given cusp $e$ is {\it rational} if it converges to
$e$, and {\it irrational} if it accumulates inside $M$.  Motivated by
problems arising from diophantine approximation, we developed in
\cite{HP} a theory of approximation of irrational geodesics by
rational ones.

As introduced in \cite{HP}, the {\it depth} $D(r)$ of a rational line
$r$ is the length of the subsegment of $r$ between the first and last
meeting point with the boundary of the maximal Margulis neighborhood
of the cusp $e$. We proved in \cite{HP} that the set of depths of
rational lines is a discrete subset of $\RR$ with finite
multiplicities. So we may define the {\it depth counting function}
$N_e:\RR\ra\NN$, with $N_e(x)$ the number of rational geodesics
whose depth is less than $x$.

Let $\widetilde{M}$ be a fixed universal cover of $M$, with covering
group $\Gamma$, and let $x_0$ be a base point in $\widetilde{M}$.
Recall that (see for example \cite{Bou1}) the {\it Poincar\'e series}
of a discrete group $G$ of isometries of $\widetilde{M}$ is
$$P(s)=\sum_{g\,\in\, G} e^{-s \,d(x_0,\,g x_0)}$$ for any $s$ in
$\RR$. This series converges if $s>\delta_G$ and diverges if
$s<\delta_G$ for some $\delta_G$ which is independant of $x_0$.
Moreover, $0<\delta_G<+\infty$ and $\delta_{G}$ is called the {\it
critical exponent} of $G$. We say that $G$ is {\it divergent} if
$P(\delta_G)$ diverges.

Choose a parabolic fixed point $\xi_0$ on the boundary
$\partial\widetilde{M}$ of $\widetilde{M}$, corresponding to $e$.
Let $\Gamma_0$ be its stabilizer. Recall that 
if $\widetilde{M}$ is a rank $1$ symmetric space (of non compact type), or if 
$\Gamma_0$ is divergent, then $\delta_G > \delta_{\Gamma_0}$ (see
\cite{DOP}. They  also give an interesting example where equality holds).

\btheo\label{theo:main_one}
If $\delta_\Gamma >\delta_{\Gamma_0}$, then
${\displaystyle \limsup_{x\ra+\infty} \frac{\log N_e(x)}{x} = \delta_\Gamma}$.
\etheo

This result on the asymptotic growth of the depths of rational lines
is easily deduced from the following statement. We define the {\it relative
  Poincar\'e series} of $(\Gamma,\Gamma_0)$ as
$$P_0(s)=\sum_{\gamma\,\in\, \Gamma_0\backslash\Gamma/\Gamma_0} e^{-s\,d(H_0,\gamma
H_0)}$$ for any $s$ in $\RR$, where $H_0$ is any fixed horosphere
centered at $\xi_0$. It is easy to see that $d(H_0,\gamma H_0)$ depends
only on the double coset of $\gamma$, and that the convergence or
divergence of the relative Poincar\'e series does not depend on the
choice of $H_0$.

\btheo\label{theo:main_two}
If $\delta_\Gamma >\delta_{\Gamma_0}$, then $P_0(s)$ converges if and 
only if $P(s)$ converges.
\etheo

These results apply in particular to the arithmetic hyperbolic
manifolds or orbifolds. For example, let $\O$ be the ring of integers
of a number field $K$, having finite group of units $O^*$ (i.e.~$K$ is
$\QQ$ or an imaginary quadratic number field $\QQ(\sqrt{-d})$). Let
$N$ be the norm on $K$, i.e.~$N(x)=x$ if $K=\QQ$ and $N(x)=|x|^2$ if
$K=\QQ(\sqrt{-d})$.  Applying Theorem \ref{theo:main_one} to the
modular orbifold $\HH^2/\PSLZ$ if $K=\QQ$ or to the Bianchi orbifold
$\HH^3/\PSLO$ if $K=\QQ(\sqrt{-d})$, one gets
 
\bcoro\label{coro:main}
Let $\varphi_{\O}(x)$ be the cardinal of the set of
$\frac{p}{q}\mbox{\rm ~mod~} \O$ with $p,q$ in $\O$,
$(p,q)=1$ and $0<N(q)\leq x$.  Then
$$\limsup_{x\ra+\infty}\frac{\log\varphi_{\O}(x)}{\log x}=2.$$ 
\ecoro

This result is well-known for $K=\QQ$, where
$\phi_{\ZZ}(x)=\sum_{k=0}^{x}\phi(k)$ and $\phi$ is the Euler
fonction. A much precise result than Corollary \ref{coro:main} is
given in the Appendix, due to K.~Belabas. Using the techniques of the
Appendix when $r_1+r_2=1$, one can hope for an analogous result for
other number fields, in connection with the counting of horoballs in
$(\HH^2)^{r_1} \times (\HH^3)^{r_2}$ under the irreducible lattice
$SL_2(\O)$ having $\QQ$-rank $1$ and $\RR$-rank $r_1+r_2$, where $r_1$
(resp.~$r_2$) is the number of real (resp.~complex) embeddings of $K$.
This will be developed later.

\medskip
{\small \flushleft{{\bf Acknowledgement.}}
We are grateful to K.~Belabas 
who provided the proof of Proposition~\ref{prop:belabas}, as well as to
E.~Fouvry and R.~de la Bret\`eche. S.~Hersonsky is grateful to the IHES 
and to the University of Orsay for their hospitality
and financial support.}

\section{Definitions}

We will use the notations and definitions of \cite{HP}, that we recall
here briefly for the sake of completeness.  We refer to \cite{HP} for
proofs and comments on these notions.

Let $M$ be a (smooth) complete Riemannian $n$-manifold with pinched
negative sectional curvature $-b^2\leq K\leq -a^2<0$. Fix a universal
cover $\widetilde{M}$ of $M$, with covering group $\Gamma$.

The boundary $\partial \widetilde{M}$ of $\widetilde{M}$ is the set of
asymptotic classes of geodesic rays in $\widetilde{M}$.  The space
$\widetilde{M}\cup \partial\widetilde{M}$ is endowed with the cone
topology. The {\it limit set} $\Lambda(\Gamma)$ is the set
$\overline{\Gamma x}\cap \partial\widetilde{M}$, for any $x$ in
$\widetilde{M}$. Let $C\Lambda(\Gamma)$ be the convex hull of the
limit set of $\Gamma$.

A point $\xi$ in $\Lambda(\Gamma)$ is a {\it conical limit point} of
$\Gamma$ if it is the endpoint of a geodesic ray in $\widetilde{M}$
which projects to a geodesic in $M$ that is recurrent in some compact
subset.  A point $\xi$ in $\Lambda(\Gamma)$ is a {\it bounded
  parabolic point} if it is fixed by some parabolic element in
$\Gamma$, and if the quotient $(\Lambda(\Gamma)-\{\xi\})/\Gamma_{\xi}$
is compact, where $\Gamma_{\xi}$ is the stabilizer of $\xi$. 

We assume that the group $\Gamma$ is {\it geometrically finite},
i.e.~that every limit point of $\Gamma$ is conical or bounded
parabolic (see \cite{Bow} for more details). We also assume that
$\Gamma$ is {\it non elementary}, i.e.~that its limit set contains at
least $3$ points.

Assume that $M$ has at least one {\it cusp} $e$, i.e.~an asymptotic
class of minimizing geodesic rays in $M$ along which the injectivity
radius goes to $0$. We say that a geodesic ray {\it converges} to $e$ if
some subray belongs to the class $e$.

Choose a parabolic fixed point $\xi_0$ on the boundary
$\partial\widetilde{M}$ of $\widetilde{M}$, which is the endpoint of a
lift of a geodesic ray converging to $e$.  Let $\Gamma_0$ be its
stabilizer in $\Gamma$. Let $H_0$ be the horosphere centered at
$\xi_0$ such that the horoball $HB_0$ bounded by $H_0$ is the maximal
horoball centered at $\xi_0$ such that the quotient of its interior by
$\Gamma_0$ embeds in $M$ under the canonical map $\widetilde{M}\ra M$.
The subset $\mbox{\rm int}(HB_0)/\Gamma_0$ of $M$ is called the
maximal Margulis neighborhood of the cusp $e$.

Since the convergence or divergence of the Poincar\'e series does not
depend on the base point $x_0$, we may assume that $x_0$ belongs to
$H_0 \cap C\Lambda(\Gamma)$.  Since the convergence or divergence of
the relative Poincar\'e series does not depend on the choice of the
horosphere, we will use this $H_0$ in the expression of $P_0(s)$ in
all that follows.

Any rational geodesic $r$ has a lift starting from $\xi_0$, which is
unique modulo the action of $\Gamma_0$. The endpoint of any such lift
is the center of an horosphere $\gamma H_0$ for some $\gamma$ in
$\Gamma$.  It follows from its definition that the depth of $r$ is
$d(H_0,\gamma H_0)$.  We proved in \cite[Lemma 2.7]{HP} that the map
$r\mapsto \Gamma_0\gamma\Gamma_0$ from the set of rational geodesics
to the set of double cosets $\Gamma_0\backslash\Gamma/\Gamma_0$ is a
bijection.  In particular the number $N_e(x)$ of rational geodesics
with depth at most $x$ is
$$N_e(x)={\rm Card}\{\Gamma_0\gamma\Gamma_0\in
\Gamma_0\backslash\Gamma/\Gamma_0\; |\; d(H_0,\gamma H_0)\leq x\}.$$

The fact that the relative Poincar\'e series $P_0(s)$ converges for
$s>\limsup_{x\ra+\infty} \frac{\log N_e(x)}{x}$ and diverges if
$s<\limsup_{x\ra+\infty} \frac{\log N_e(x)}{x}$ is then easily seen.
In particular, Theorem \ref{theo:main_one} follows from Theorem
\ref{theo:main_two}.

\medskip
Let $\O$ be the ring of integers of a number field $K$, where $K$ is
either $\QQ$ or an imaginary quadratic number field $\QQ(\sqrt{-d})$,
with $d$ a positive square free integer. We use the upper half-space
models for the real hyperbolic spaces. Consider the cusp $e$ in the
orbifolds $\HH^2/\PSLZ$ if $K=\QQ$, or $\HH^3/\PSLO$ otherwise,
corresponding to the parabolic fixed point $+\infty$. We proved in
\cite[Section 2.3]{HP} (with the obvious adaptation to the case of orbifolds) that
the rational lines $r$ are in one-to-one correspondance with the
fractions $\frac{p}{q}$ modulo the additive group $\O$, with the depth
of $r$ being $\log |q|^2$, if this fraction is written with relative
prime numerator and denominator.  Hence Corollary \ref{coro:main}
follows from Theorem \ref{theo:main_one}.

\section{Proofs}

With the notations and assumptions of the previous section,
we start with two lemmae.

\blemm\label{lemm:close_to_pick} 
There exists a constant $C_1\geq 0$, such that every double coset
$\Gamma_{0}\gamma\Gamma_{0}$ in
$\Gamma_{0}\backslash\Gamma/\Gamma_{0}$ has a representative $\gamma$
which satisfies
$$|d(H_0,\gamma H_0)-d(x_0,\gamma x_0) |\leq C_1.$$
\elemm

\dem 
Choose the identity as the representative of the trivial double coset.
Let $\gamma$ be in $\Gamma-\Gamma_0$.  Let $p_0$ in $H_0$ and $p_1$ in $\gamma
H_0$ be such that the segment $[p_0,p_1]$ is the (unique) common
perpendicular to $H_0$ and $\gamma H_0$. In particular, $d(H_0,\gamma
H_0)=d(p_0, p_1)$, and $p_0,p_1$ lie on the geodesic line between the
centers of $H_0$ and $\gamma H_0$, so that $p_0,p_1$ both belong to
the $\Gamma_0$-invariant subset
$H_0 \cap C\Lambda(\Gamma)$. Since $\xi_0$ is a bounded parabolic fixed point,
the quotient $(H_0\cap C\Lambda(\Gamma))/\Gamma_0$ is compact, hence has
diameter bounded by $C'_1\geq 0$.

Since $x_0$ belongs to $H_0 \cap C\Lambda(\Gamma)$, there exists
$\alpha$ in $\Gamma_0$ so that $d(p_0, \alpha x_0)\leq C'_1$, 
and $\beta$ in $\Gamma_0$ so that $d(\gamma^{-1}p_1, \beta x_0)\leq C'_1$.

Since $ \alpha x_0$ lies on $H_0$ and $\gamma\beta x_0$ on $\gamma
H_0$, we have
$$d(H_0,\gamma
H_0)=d(p_0, p_1)\leq d( \alpha x_0, \gamma\beta x_0).$$
Conversely, by the triangular inequality,
$$d( \alpha x_0, \gamma\beta x_0)\leq d( \alpha x_0,p_0) + d(p_0,p_1)
+ d(p_1, \gamma\beta x_0)\leq 2C'_1+d(H_0,\gamma H_0).$$ Hence the
representative $\alpha^{-1}\gamma\beta$ of the double coset
$\Gamma_0\gamma\Gamma_0$ satisfies the condition of the Lemma with
$C_1 = 2C'_1$.  \eop{\ref{lemm:close_to_pick}}

\bigskip
Note that by discreteness, the set of representatives as in the Lemma
of a given double coset is finite.  From now on, we will denote by the
same letter a double coset and such a representative.

The following proposition is well-known (see for example the proof of
\cite[Lemme 4]{DOP}).

\blemm\label{lemm:thin_triangles}
There exists a constant $C_2\geq 0$ (depending only on the upperbound
on the curvature of $M$) such that the following holds.  Let $H,H'$ be
horospheres in $\widetilde{M}$ bounding disjoint horoballs.  Let
$[p,p']$ be the common perpendicular segment, with $p$ in $H$ and $p'$
in $H'$.  For every $x$ in $H$ and $x'$ in $H'$, 
$$|d(x,x')-(d(x,p) + d(p,p') +d(x,p))|\leq C_2.$$
\elemm

\addtocounter{fig}{1}
\[
\begin{array}{c}
\mbox{\begin{picture}(0,0)%
\includegraphics{thin.pstex}%
\end{picture}%
\setlength{\unitlength}{3158sp}%
\begingroup\makeatletter\ifx\SetFigFont\undefined
\def\x#1#2#3#4#5#6#7\relax{\def\x{#1#2#3#4#5#6}}%
\expandafter\x\fmtname xxxxxx\relax \def\y{splain}%
\ifx\x\y   
\gdef\SetFigFont#1#2#3{%
  \ifnum #1<17\tiny\else \ifnum #1<20\small\else
  \ifnum #1<24\normalsize\else \ifnum #1<29\large\else
  \ifnum #1<34\Large\else \ifnum #1<41\LARGE\else
     \huge\fi\fi\fi\fi\fi\fi
  \csname #3\endcsname}%
\else
\gdef\SetFigFont#1#2#3{\begingroup
  \count@#1\relax \ifnum 25<\count@\count@25\fi
  \def\x{\endgroup\@setsize\SetFigFont{#2pt}}%
  \expandafter\x
    \csname \romannumeral\the\count@ pt\expandafter\endcsname
    \csname @\romannumeral\the\count@ pt\endcsname
  \csname #3\endcsname}%
\fi
\fi\endgroup
\begin{picture}(3450,3450)(676,-2611)
\put(2476,-1411){\makebox(0,0)[lb]{\smash{\SetFigFont{11}{13.2}{rm}$p$}}}
\put(1126,-1261){\makebox(0,0)[lb]{\smash{\SetFigFont{11}{13.2}{rm}$x$}}}
\put(2851,314){\makebox(0,0)[lb]{\smash{\SetFigFont{11}{13.2}{rm}$x'$}}}
\put(3376,-136){\makebox(0,0)[lb]{\smash{\SetFigFont{11}{13.2}{rm}$p'$}}}
\end{picture}
}\\
\\ 
\hbox{\rm Figure \arabic{fig} ~: The quasi-geodesic.}
\end{array}
\]

\dem  Up to replacing $H,H'$ by inside concentric horospheres at
distance $1$, we may assume that $d(H,H')\geq 1$. By the convexity of
horoballs, the piecewise geodesic $[x,p]\cup[p,p']\cup[p',x']$ has
angles at least $\frac{\pi}{2}$ at $p$ and at $p'$. Thus, since
$d(p,p')\geq 1$, it is a quasi-geodesic, and the result follows for
example from \cite[Chapter 3]{GH}. 
\eop{\ref{lemm:thin_triangles}}

\bigskip
\noindent{\bf Proof of Theorem~\ref{theo:main_two}} If $f,g$ are maps from
an interval $I$ in $\RR$ to $\RR\cup\{+\infty\}$, write $f\asymp g$ if
there exists a finite constant $c>0$ such that $\frac{1}{c}g(s)\leq f(s)\leq
c g(s)$ for all $s$ in $I$.  We write the Poincar\'e series as
follows.
$$P(s)=\sum_{\gamma\,\in\,\Gamma_0\backslash\Gamma/\Gamma_0}\;
\sum_{\alpha,\,\beta\in\Gamma_0} e^{-s\,d(x_0,\alpha\gamma\beta x_0)}.$$

Note that $d(x_0,\alpha\gamma\beta x_0)= d(\alpha^{-1}x_0,\gamma\beta
x_0)$. The representatives $\gamma$ of the (non trivial) double cosets have
been chosen so that $\gamma x_0$ lies at distance less than a constant
$C'_1$ from the endpoint on $\gamma H_0$ of the common perpendicular
segment between $H_0$ and $\gamma H_0$. Applying
Lemma~\ref{lemm:close_to_pick} and Lemma~\ref{lemm:thin_triangles}, 
we obtain
$$P(s)\asymp \sum_{\gamma\,\in \,\Gamma_0\backslash\Gamma/\Gamma_0}\;
\sum_{\alpha,\,\beta \in \Gamma_{0}} e^{-s\left(d(\alpha^{-1}x_0,x_0) +
d(H_0,\gamma H_0) + d(\gamma x_0,\gamma\beta x_0)\right)}$$
$$=\left(\sum_{\gamma\,\in \,\Gamma_0\backslash\Gamma/\Gamma_0} 
e^{-s\,d(H_0,\gamma H_0)}\right)\left(
\sum_{\alpha\in \Gamma_{0}} e^{-s\,d(x_0,\alpha x_0)}\right)^2.$$

If $s>\delta_{\Gamma_0}$, then the Poincar\'e series of $\Gamma_0$
converges at $s$. Hence $P\asymp P_0$ on the interval
$]\delta_{\Gamma_0},+\infty[$. Theorem \ref{theo:main_two} follows.
\eop{\ref{theo:main_two}}

\section{Appendix, by K.~Belabas}
 
For all undefined objects and unproved results in what follows, see
\cite{Nar} for instance.  

Let $K$ be a number field, $N=N_K$ the norm on $K$, $\zeta_K$ the
Dedekind zeta function and Res$_K={\rm Res}(\zeta_K,s=1)$ the residue
of $\zeta_K$ at the unity, $\O=\O_K$ the ring of integers of $K$,
$\O^*$ the units of $\O$, $h_K$ the class number of $K$, $W_K$ the
number of roots of unity in $K$, $R_K$ the regulator of $K$, and $D_K$
the discriminant of $K$. Let $\mu$ be the M\"obius function, defined
by $\mu(\O)=1$, $\mu(I)=(\-1)^k$ if the ideal $I$ in $\O$ is the
product of $k$ distinct prime ideals, and $\mu(I)=0$ if $I$ is
divisible by the square of a prime ideal. It satisfies
$\sum_{I\,|\,J}\mu(I)= 1$ if $J=\O$, and $0$ otherwise.  Assume that
$\O^*$ is finite ($K$ is $\QQ$ or an imaginary quadratic field), so
that Card$(\O^*)=W_K$ and $R_K=1$.

For $n$ a positive integer, define
$$\phi(x)=\phi_\O(x)={\rm ~Card}\{\frac{p}{q} {\rm ~mod~} \O,
(p,q)\in\O\times(\O-\{0\}), (p,q)=1 {\rm ~and~} N(q)\leq x\}.$$

The following proposition was explained to us by K.~Belabas, who said
it might be already known. For the sake of completeness, we include
its proof here.

\bprop\label{prop:belabas}
There exists $\epsilon>0$ such that
$$\phi(x)=\frac{{\rm Res}_K}{2\,h_K\,\zeta_K(2)} \;x^2 + 
O(x^{2-\epsilon}).$$
\eprop

\noindent{\bf Proof} (K.~Belabas).
Since two irreducible fractions are equal if and only if the numerators
and denominators are multiplied by the same unit, one has
$$\phi(x)=\frac{1}{W_K}{\rm ~Card}\{(q,p{\rm ~mod~}(q)), q\neq 0,
(p,q)=1 {\rm ~and~} N(q)\leq x\}.$$
Hence $W_K\phi(x)$ is equal to
$$
\sum_{(q,p{\rm ~mod~}(q)), \,q\neq 0, \,N(q)\leq x, \,(p,q)=1} 1 = 
\sum_{(q,p{\rm ~mod~}(q)), \,q\neq 0, \,N(q)\leq x} \;\sum_{I\,|\,(p,q)}
\mu(I) = \sum_{I} \mu(I)f(I)$$ 
where $I$ ranges over the ideals of $\O$ and
$$f(I)=\sum_{(q,p{\rm ~mod~}(q)), \,q\neq 0,
  \,I\,|\,(p,q), \,N(q)\leq x} 1= \sum_{q\in I,
  \,q\neq 0,\,N(q)\leq x} \;\sum_{p\in I/(q)} 1= W_K
\sum_{(q)\subset I, \,N(q)\leq x} \frac{N(q)}{N(I)}$$
since $N(J)={\rm Card~}\O/J$, where $(q)$ is a non zero principal
ideal, and since a generator of a principal ideal is uniquely defined up
to units.

\blemm\label{lem:belabasun}
If $S(x)={\rm Card} \{ (q)\subset I, N(q)\leq x\}$, then 
$S(x)=\frac{{\rm Res}_K}{h_KN(I)}x+ O((\frac{x}{N(I)})^{1-\epsilon})$.
\elemm

\dem Note that $(q)\subset I$ if and only if $(q)=IJ$ for some ideal
$J$ in $\O$, and $N(IJ)=N(I)N(J)$. Hence $$S(x)=\sum_{J\in[I]^{-1},\,
N(J)\leq \frac{x}{N(I)}} 1$$ where $[I]^{-1}$ is the inverse of the
class of $I$ in the class group. With $\epsilon=1/[K:\QQ]$, the result
then follows from \cite[theo.~7.6 page 361]{Nar} (for example) for the
main term and from \cite{Tat} for the error term.
\eop{\ref{lem:belabasun}}

\blemm\label{lem:belabasdeux}
If $T(x)=\sum_{(q)\subset I, \,N(q)\leq x} N(q)$, then 
$T(x)=\frac{{\rm Res}_K }{2h_KN(I)}x^2+ 
O(\frac{x^{2-\epsilon}}{N(I)^{1-\epsilon}})$.
\elemm

\dem 
This is immediate by applying Fubini and using the previous lemma.
\eop{\ref{lem:belabasdeux}}

\bigskip
Now, since $\zeta_K(2-\epsilon)$ converges for $\epsilon<1$,
and since by \cite[page 326]{Nar},
$$\sum_{I} \frac{\mu(I)}{N^s(I)} =\frac{1}{\zeta_K(s)},$$ one gets
from Lemma \ref{lem:belabasdeux}
$$\phi(x)=
\sum_{I} \mu(I)\frac{{\rm Res}_K \;x^2}{2 h_K N^2(I)}
+O(x^{2-\epsilon})= 
\frac{{\rm Res}_K}{2\,h_K\,\zeta_K(2)} \;x^2 + O(x^{2-\epsilon}).$$
\eop{\ref{prop:belabas}}

Since Res$_K=\frac{2^{r_1}(2\pi)^{r_2}h_KR_K}{W_k\sqrt{|D_K|}}$, one
has for $K=\QQ(\sqrt{-d})$, with $D=d$ if $d\equiv 3 {\rm ~mod~} 4$ and
$D=4d$ otherwise, with $w=4,6$ if $d=1,3$ and $w=2$ otherwise,
$$\phi_{\O_{-d}}(x)=
\frac{\pi}{w\,\zeta_{\QQ(\sqrt{-d})}(2)\sqrt{D}}\;x^2 + o(x^2).$$
For $K=\QQ$, the formula gives the well-known result (see \cite{Nar} 
for instance)
$$\phi_\ZZ(x)=\frac{3}{\pi^2}x^2 + o(x^2).$$

\bigskip
\noindent {\small
\begin{tabular}{l} 
Department of Mathematics\\ 
Ben Gurion University\\ BEER-SHEVA, Israel\\
{\it e-mail: saarh@math.bgu.ac.il}
\end{tabular}
\\
 \mbox{}
\\
 \mbox{}
\\
\begin{tabular}{l} Laboratoire de Math\'ematiques UMR 8628 CNRS\\
Equipe de Topologie et Dynamique (B\^at. 425)\\
Universit\'e Paris-Sud \\
91405 ORSAY Cedex, FRANCE.\\
{\it e-mail: Frederic.Paulin@math.u-psud.fr}
\end{tabular}
}

\end{document}